\def\AA         {{\bf A}}
\def\ZZ         {{\bf Z}}
\def\RR         {{\bf R}}
\def\CC         {{\bf C}}

\def\GG         {{\bf G}}
\def\LL         {{\bf L}}
\def\QQ         {{\bf Q}}
\def\PP         {{\bf P}}
\def\A          {{\cal A}}
\def\L          {{\cal L}} 
\def\M          {{\cal M}} 
\def\S          {{\cal S}}
\def\P          {{\cal P}}
\def\F          {{\cal F}}
\def\D          {{\cal D}} 
\def\O          {{\cal O}}
\def\V          {{\cal V}}

\def\dim        {{\rm dim}}

\def\Spec       {{\rm Spec}}
\def\Char       {{\rm Char}}
\def\log        {{\rm log}}
\def\Hom        {{\rm Hom}}

\documentstyle[twoside,12pt]{article}
\newtheorem{prop}{Proposition}[section]
\newtheorem{dfn}[prop]{Definition}   
\newtheorem{theo}[prop]{Theorem}
\newtheorem{conj}[prop]{Conjecture}
 
\newtheorem{coro}[prop]{Corollary} 
\newtheorem{lem}[prop]{Lemma}   
\newtheorem{exam}[prop]{Example}
\newtheorem{prob}[prop]{Problem}

\title{Homotopy groups of complements to ample divisors}
\author{A.Libgober\\
\small Department of Mathematics \\
\small University of Illinois at Chicago\\
\small 851 S.Morgan Str. Chicago, Illinois, 60607 \\
\small e-mail: libgober@math.uic.edu \\}

\begin{document}

\date{}

\maketitle

\begin{abstract}{We study the homotopy groups of complements to 
reducible divisors on non-singular projective varieties
with ample components and isolated non normal crossings. 
We prove a vanishing 
theorem generalizing conditions for commutativity of the fundamental 
groups. We also discuss the calculation of supports of  
non vansihing homotopy groups as modules over the fundmanetal group
in terms of the geometry of the locus of non-normal crossings.
We review previous work on the
local study of isolated non-normal crossings and relate 
the motivic zeta function to the local polytopes of quasiadjunction. 
As an application, we obtain information about the support 
loci of homotopy groups of arrangements of hyperplanes.}
\end{abstract}

\section{Introduction} An interesting problem in the study of the topology 
of algebraic varieties is understanding 
the fundamental group of the complement to a divisor 
on a non-singular algebraic variety in terms of the 
geometry of the divisor.
Works of Abhyankar (\cite{Abh}) and Nori (\cite{Nori}) show 
that, if $C$ is an irreducible curve on a 
non singular algebraic surface $X$, then for some effective constant
$F(C)$ depending on the local type of singularities of $C$, the 
inequality $C^2>F(C)$ implies that the 
kernel $\pi_1(X-C) \rightarrow \pi_1(X)$ is a central extension.
For example, if $X$ is simply connected, then $\pi_1(X-C)$ is 
abelian. Historically, such results were originated in the so-called 
Zariski problem and we refer to \cite{Fulton} for a survey.
The case of non-abelian fundamental groups of complements, 
notably when $X={\bf P}^2$, is also very interesting.
The geometric information, such as the dimensions of the linear 
systems defined by singularities 
of the curve, becomes essential in descriptions of  
fundamental groups and their invariants 
(cf. \cite{mina}, \cite{arcatame2} \cite{abcov}). 
Recently, analogous questions about fundamental groups 
of the complements in the case when $X$ is 
symplectic began to  attract attention as well
(cf. \cite{symplectic}).

In the present work, we shall show that, in appropriate settings, 
the relationship between the topology of the complement and the 
geometry of the divisor can be extended to some higher homotopy groups.
Some work in this direction already was done. In \cite{Annals}, we show that 
if $V$ is a hypersurface in ${\bf C}^{n+1}$ with isolated singularities 
which compactification in ${\bf P}^{n+1}$ is transversal to 
the hyperplane at infinity, then the first homotopy groups 
of the complement are the following: 
\begin{equation}\label{pi1P}
\pi_1(\CC^{n+1}-V)=\ZZ,  \ \ \ \pi_i(\CC^{n+1}-V) \ {\rm for} \ 2 \le i \le n-1
\end{equation}
Moreover, the next homotopy group $\pi_n(\CC^{n+1}-V)$
depends on the local type 
of the singularities and also on the geometry of a collection 
of singularities as a finite subset in $\CC^{n+1}$.
It can also be described via a generalization of the van Kampen 
procedure in terms of pencils of hyperplane sections (cf. \cite{Annals},
\cite{denis}).
Recently, homotopy groups of arrangements were considered 
in \cite{DP} and  \cite{SP}.

Below, we shall extend these results in two directions.
On the one hand, we shall consider complements on arbitrary algebraic 
varieties rather than just in projective space. The latter case, 
however, appears to be the most important one due to a variety 
of interfaces with 
other areas -- e.g. the study of arrangements of hyperplanes.
On the other hand, we do not assume here that $V$ has isolated 
singularities, but rather that the divisor $D$ has normal crossings 
except for finitely many points. The effect of this is that 
the fundamental group, which plays the key role in the description 
of higher homotopy,  may be abelian rather than cyclic, as 
is the case in (\ref{pi1P}),
and the theory which we obtain is {\it abelian} rather than {\it cyclic}.

In the next section, we prove the triviality of   
the action of the fundamental group on higher homotopy groups
in certain situations (cf. theorem \ref{triviality}).
This implies that all information about homotopy groups 
in these cases is homological ($\pi_1$ in this situation 
is automatically abelian). In some instances, as result of homological
calculations, one obtains a {\it vanishing}
of homotopy groups in certain range.
In particular if $D$ is a divisor in $\PP^{n+1}$ having only 
isolated non normal crossings and the number of components 
greater than $n+1$ then $\pi_i(\PP^{n+1}-D)=0$ 
in the range $2 \le i \le n-1$. The results of section 
\ref{action} also isolate first non-trivial
homotopy group in the sense that it is a non-trivial 
$\pi_1$-module.

 In section \ref{charvarsection} we define our main invariant
of the homotopy group, i.e. a  sub-variety of the 
spectrum of the group ring of the fundamental group which is the support
of the first non trivial homotopy group considered as the module over $\pi_1$. 
We call these sub-varieties {\it characteristic}
and show that these sub-varieties are related 
to the jumping loci for the cohomology of local systems. The latter 
have a very restricted structure (cf. \cite{Arapura}),-- i.e. they 
are unions of translated by points of finite order subgroups 
which also suggest the numerical data 
that describe these varieties completely.

The methods of obtaining numerical data specifying 
the characteristic varieties from the geometry of the divisor 
are discussed in the section \ref{conjectures}. This is done
by using the Hodge theory of abelian covers, which is studied in 
section  \ref{covers}, and relies on our local study of the isolated 
non normal crossings in \cite{innc} and \cite{loctop}. 
The results of these papers are discussed in the section \ref{review}
where, among other things, we compare cyclic theory of 
isolated singularities with abelian theory of isolated-non normal crossings 
having more than one components. 
In section \ref{examples}, we review cases to which the results 
of section \ref{conjectures} can be applied. In particular, 
the Kummer configuration yields an arrangement of planes in $\PP^3$ 
with non-trivial $\pi_2$ of the complement which we calculate.
In the final section, we show the relationship between the 
invariants of the homotopy groups and the motivic zeta function 
of Denef-Loeser.

Part of this work was done during my visit to University of Bordeaux
to which I wish to express my gratitude. I particularly want 
to thank Alex Dimca and Pierrette Cassou-Nogues for their
hospitality. I also want express my gratitude to the organizers of the 
Sapporo meeting on Singularities where parts of the results 
of this paper were presented.

\section{Action of the fundamental group on homotopy groups.}
\label{action}

In this section we discuss homotopy groups, 
in a certain range of dimensions, for  
a class of quasi-projective varieties. This is done in 
two steps. Firstly, we show that these varieties
support a trivial action of $\pi_1$ (in particular are {\it nilpotent}
in certain range). Secondly, we use homological calculations to 
determine these homotopy groups and to describe cases when 
homotopy groups vanish. 

Recall that homotopy groups $\pi_n(X,x)$ of a topological space $X$ 
are $\pi_1(X,x)$-modules, with the action given by 
the ``change of the base point'' (cf. \cite{Spanier}). 
In the case when $\pi_i(X)=0$ for $1 <i < n$, this action on 
$\pi_n(X)$, which is isomorphic to $H_n(\tilde X)$ where $\tilde X$ 
is the universal cover of $X$, coincides with the action 
of the fundamental group on the homology 
of the universal cover via deck transformations.
A topological space is called $k$-simple 
if the action of the fundamental group on $\pi_i(X)$ is trivial 
for $i \le k$ (cf. \cite{Spanier}).

Examples of $k$-simple spaces appearing 
naturally in algebraic geometry 
are the following. Locally, they come up when 
one looks at the complement 
to a union of germs of divisors in $\CC^{n+1}$ 
forming an isolated non-normal crossing. This situation 
was studied in \cite{innc}.
More generally (cf. \cite{loctop}), 
instead of divisors in $\CC^{n+1}$, one can look at the union of 
germs  of divisors in a germ of a complex space $Y$ having a link 
which is $(\dim Y-2)$-connected. Examples of such local singularities 
are provided by the cones over a normal crossings divisor in $\PP^{n}_{\CC}$, 
in particular by cones over generic arrangements of hyperplanes.
These $k$-simple spaces  ($k=\dim Y-2$) are, of course, Stein spaces and their
theory will be reviewed in 
section \ref{review}. Global $k$-simple examples 
are given by the following quasi-projective varieties: 

\begin{theo}\label{triviality} Let $X$ be a {\it simply connected}
projective manifold and $D=\bigcup D_i$ 
be a divisor {\it with normal crossings} 
such that its all components $D_i$ are smooth and ample.
Then $\pi_1(X-D)$ is abelian and its action on $\pi_i(X-D)$ is trivial 
for $2 \le i \le {\rm dim}X-1$. 
\end{theo}

The proof is similar to the one presented in the local case in \cite{innc}.
It uses the reduction to the case of normal crossings divisors
using Lefschetz hyperplane section theorem and then 
surjectivity of $\pi_i(D_i-\bigcup_{j \ne i} D_j) \rightarrow \pi_i(X-D)$
which follows from ampleness of the components $D_i$.

This theorem reduces the calculation of the homotopy group to
the calculation of {\it homology} of the complements. The latter
can be done using the  exact sequence:

\begin{equation}\label{homologycomplement}
H_2(X) \rightarrow H_2(X,X-D) \rightarrow H_1(X-D) 
\rightarrow H_1(X) \rightarrow H_1(X,X-D)
\end{equation}
and the isomorphisms $$H_j(X,X-D)=H^{2n+2-j}(D)$$
We obtain hence:
\begin{coro}\label{h1}
 Let $H={\bf Z}^N$ be a free abelian 
group generated by components of the divisor $D$. 
Let $$h: H_2(X,{\bf Z}) \rightarrow H$$ given by $a \rightarrow 
\sum (a,D_i)D_i$ where $a \in H_2(X), D_i \in H^2(X)$ and 
$(a,D)$ is the Kronecker pairing. Then $\pi_1(X-D)={\rm Coker} h$.
For example, if $X={\bf P}^{n+1}$ and one of the 
components $D_i$ ($i=1,..,r+1)$ 
is a hyperplane, then $\pi_1(X-D)={\bf Z}^{r}$. 
Let $X$ be a hypersurface in ${\bf P}^{n+1}$ and $D$ be 
a union of $r+1$-hyperplanes. Then $H_1(X-D)=\ZZ^r$.

\end{coro}
The following result can be used for the calculation of 
the homology of some branched covers of $X$:

\begin{coro} Let $D_i \in \vert \L_i^{m_i} \vert (i=1,...,r)$ such that 
$D_i$ are divisors on $X$ having isolated non-normal crossings 
and $D_i$ is the zero set of $f_i \in H^0(X,\O(D_i))$. Let $s_i
 \in H^0(X,\L_i) \ne 0$. Then $\bar U_{m_1,..,m_r}$ given in 
the total space of $\oplus \L_i$ 
 by $s_i^{m_i}=f_i$ is  
the cover corresponding the surjection:
$\phi: H_1(X-D) \rightarrow G=\oplus_{i=1} \ZZ/m_i\ZZ$ 
The projection $\bar U_{m_1,..,m_r} \rightarrow X$ induces the 
isomorphism: $H_i(\bar U_{m_1,...,m_r}) \rightarrow H_i(X)$ 
for $i \le n-1$.

\end{coro} 
The next theorem is an abelian version of 
the result in \cite{Annals} and 
identifies ``the first non-trivial 
homotopy group'' in the sense of \cite{Annals}. 

\begin{theo}\label{arrangements} 
 $\rm (a)$ Let $X={\bf P}^{n+1}$ and $D$ be an arrangement of 
$r+1$ hypersurfaces as in corollary (\ref{h1})
(i.e., such that one of the hypersurfaces has a degree 1) 
and having finitely many non-normal crossings. 
Then $\pi_i({\bf P}^{n+1}-D)=0$ for $2 \le i \le n-1$. 
If all intersections are the normal crossings, then the 
$\pi_n(\PP^{n+1}-D)=0$.
\smallskip 
\par \noindent $\rm (b)$ Let $V$ be a complete intersection in $\PP^N$ and 
$\dim V=n+1$. Let $D$ be the arrangement of $r+1$ hyperplane sections
of having isolated non-normal crossings. 
Then $\pi_1(V-D)=\ZZ^r$ and 
$\pi_i(V-D)=0$ for $2 \le i \le n-1$.
 \end{theo} 

\noindent {\it Proof.} Consider first $\rm (a)$. 
The claimed vanishing 
is a consequence of 
the Lefschetz hyperplane section theorem
(cf. \cite{GM}) 
and the second part of $\rm (a)$. The first part  
follows by induction, with the inductive step 
being the vanishing of $\pi_n(\PP^{n+1}-D)$ where $D$ is 
an arrangement of hypersurfaces with normal crossings.
Taking into account the triviality of the 
action of $\pi_1(\PP^{n+1}-D)$ on $\pi_n$, the claim is a 
consequence of the exact sequence (cf. \cite{Brown}):

\begin{equation}
H_{n+1} (\PP^{n+1}-D) \rightarrow  
H_{n+1} ({\bf Z}^{r}) \rightarrow 
\pi_{n}(\PP^{n+1}-D)_{{\bf Z}^{r}} \rightarrow
\end{equation}
$$ \rightarrow H_{n}(\PP^{n+1}-D)\rightarrow 
H_{n} ({\bf Z}^{r}) \rightarrow 0$$   
and the calculation of the homology of $\PP^{n+1}-D$.
The latter can be done using Mayer Vietors spectral sequence
(cf. \cite{innc}). The proof of $\rm (b)$ is similar.

\section{Characteristic varieties of homotopy groups}
\label{charvarsection}

In this section, we study the support of the first homotopy 
group of quasi-projective varieties from section \ref{action} on which 
the action of $\pi_1$ fails to be trivial. This support is 
a subvariety of $\Spec \CC[\pi_1]$, which we call {\it the 
characteristic variety}.
We show that 
in the range $2 \le  i \le k-1$, in which the action of $\pi_1$ 
on $\pi_i$ is trivial, 
the homology $H_i$ of 
the local systems, corresponding to the points of the algebraic 
group $\Spec \CC[\pi_1]$ different from the identity, 
is trivial. 
Moreover, the first homotopy group outside this range, 
i.e. $\pi_k$, determines 
the homology $H_k$ of the local systems.
Vice versa, the (co)homology of local systems determines the 
support of $\pi_k \otimes \CC$ as $\pi_1$-module. 
This yields, in the algebro-geometric context, a ``linear'' structure
of the characteristic varieties.

\begin{theo} \label{hurewicz}
Let $X$ be a topological space with an abelian fundamental 
group $\pi_1(X)=A$. Assume that for an ideal $\wp$ in ${\bf C}[A]$
the localization of the homotopy groups is trivial for $2 \le i <k$.
$\pi_i(X)_{\wp}=0$. Then $H_i(\tilde X)_{\wp}=0$ for $1 \le i <k$ 
and $H_k(\tilde X)_{\wp}=\pi_k(X)_{\wp}$ 

\end{theo} 

\noindent {\it Sketch of the proof}. The universal cover $\tilde X$ of $X$
is a simply connected space on which $A$ acts freely. For such 
a space, the group $A$ acts on $H_j(\tilde X,{\bf C})$ for any $j$ and 
on the homotopy groups $\pi_j(\tilde X,\tilde x_0)=\pi_j(X,x_0)$ 
($j \ge 2$)
so that the Hurewicz map: $\pi_j(\tilde X) \rightarrow H_j(\tilde X)$
is $\pi_1(X)$-equivariant (cf. \cite{Spanier} Ch.7, Cor. 3.7). 

Let us consider a simply connected 
CW-complex $Y$ on which an abelian group $A$ acts freely. 
The group $A$ then acts on the homotopy groups via composition
of the map $\pi_n(Y,x) \rightarrow \pi_n(Y,a(x))$ and the identification
$\pi_n(Y,a(x))$ and $\pi_n(Y,x)$, which is independent 
of the choice of a path connecting $x$ and $a(x)$ due 
to $\pi_1(Y)=0$. The claim is that, if $\pi_i(Y)_{\wp}=0$ for $1 <i \le n-1$,
then 
$\pi_n(Y)_{\wp}=H_n(Y)_{\wp}$
 where $\wp \in {\rm Spec}\CC [A]$ is such that $\wp$ is not in 
the maximal ideal of the identity of group ${\rm Spec}\CC [A]$. 
The theorem above will follow for $Y=\tilde X$ and $G=\pi_1(X)$.

Consider the fibration of path space $Maps(I,Y) \rightarrow Y \times Y$.
This fibration is equivariant (where the action on $Y \times Y$ 
is diagonal). Space $Map(I,Y)$ is homotopy equivalent to $Y$. 
We have the spectral sequence: 
 $$E_2^{p,q}: H_p(Y \times Y,H_q(\Omega Y)) \rightarrow H_{p+q}(Y)$$ 
This spectral sequence is equivariant with the action on the fiber 
given by $gv=p_*g_*(v)$ where 
$g_*: H_i(\Omega_xY) \rightarrow H_i(\Omega_{gx}Y)$ and $p_*$ 
is the Gauss-Manin connection identifying homology 
of different fibers of fibration.
Localizing at $\wp$, due to $(n-1)$ simplicity of $Y$, we
obtain that the terms with $0< p \le n-1$ and $0<q \le n-2$ are zeros. 
In localized spectral sequence we can identify the map 
$H_n(Y)_{\wp} \rightarrow E_{\infty}^{n,0}=
Ker d_n^{n,0}: H_n(Y \times Y)_{\wp} 
\rightarrow H_{n-1}{\Omega Y}_{\wp}$ with the map 
$i_{\Delta}: H_n(Y)_{\wp} \rightarrow 
H_n(Y \times Y)_{\wp}$ corresponding to a diagonal embedding.
Moreover, $d_n^{n,0}$ is surjective 
(since $H_{n-1}(Y)_{\wp}=0$). Hence we have an exact sequence:
$$0 \rightarrow Im(i_{\Delta})_{\wp} \rightarrow H_n(Y \times Y)_{\wp}
\rightarrow H_{n-1}(\Omega Y)_{\wp} \rightarrow 0$$
and therefore $H_n(Y)_{\wp}=H_{n-1}(\Omega Y)_{\wp}=\pi_n(Y)_{\wp}$.

We shall apply this theorem to $(n-1)$-simple spaces. 
For such a space the support of $\pi_i(X) \otimes_{\ZZ} \CC$ 
as a $\CC[\pi_1(X)]$ 
module belongs for $2 \le i \le n-1$ to the maximal ideal 
of the identity of the 
group $\Spec \CC[\pi_1(X)]=\Char [\pi_1(X)]$. This maximal 
ideal is just the augmentation ideal of the group ring.
Hence the localization at a prime ideal not belonging to 
the maximal ideal of the identity satisfies (after tensoring with $\CC$) 
the assumption of the theorem \ref{hurewicz}. This allows, for 
$(n-1)$-simple spaces, to express the homology of the local systems
in terms of the homotopy groups $\pi_n(X)$:

\begin{theo} \label{locsys}
 Let $\rho \in {\rm Char}\pi_1(X)$ be a 
character of the fundamental group different from the identity and 
let ${\CC}_{\rho}$ be $\CC$ considered as $\CC [\pi_1(X)]$ module 
via the character $\rho$. Then 
$$H_i(X,\rho)=0 \ (i \le n-1) \ \ \ H_n(X,\rho)=\pi_n(X) \otimes_{\CC [\pi_1(X)]} {\CC}_{\rho}$$
\end{theo} 

\noindent {\it Proof.} The proof is similar to the one in the case 
when $X$ is a complement to a plane curve (cf. \cite{abcov}) and 
the local case (cf. \cite{innc}). Consider the spectral sequence
(cf. \cite{cartaneilenberg}, ch. XVI th.8.4):
  $$H_p(\pi_1(X),H_q(\tilde X)_{\rho}) \Rightarrow H_{p+q}(X,\rho)$$
where $H_*(\tilde X)_{\rho}$ is the homology of the 
complex $C(\tilde X) \otimes_{\ZZ} {\CC}$ with the 
action of $\pi_1(X)$ given by $g(e \otimes \alpha)=g \cdot e \otimes 
\rho(g^{-1})\alpha$. We can localize this spectral sequence 
at the maximal ideal $\wp_{\rho}$ of ${\rm Spec}[\pi_1(X)]$ corresponding 
to the character $\rho$. The resulting spectral sequence has 
$E_2^{i,j}=0$ for $1 \le j \le n-1$. The exact sequence of low degree
terms yields:  
$H_n(X,\rho)=H_n(\tilde X) \otimes_{\CC [\pi_1(X)]} {\CC}_{\rho} $
which together with the theorem \ref{hurewicz} proves the claim.

Now we are ready to define the main invariant.

\begin{dfn} The $k$-th characteristic variety $V_k(\pi_n(X))$
of the homotopy group $\pi_n(X)$
is the zero set of the $k$-th Fitting ideal of $\pi_n(X)$, i.e. 
the zero set of minors of order $(n-k+1) \times (n-k+1)$ of $\Phi$ in a 
presentation 
$$\Phi: \CC [\pi_1(X)]^m \rightarrow \CC [\pi_1(X)]^n 
\rightarrow \pi_n(X) \rightarrow 0$$
 of $\pi_1(X)$ module $\pi_n(X)$ via generators and relations. 
Alternatively (cf. theorem \ref{locsys})
outside of $\rho=1$, $V_k(\pi_n(X))$ is the set 
of characters $\rho \in \Char [\pi_1(X)]$ such that  
$\dim H_n(X,\rho) \ge k$.
\end{dfn}



Theorem \ref{locsys} combined with the results of \cite{Arapura}
yields the following strong structure property (for possibly 
non-essential characters):

\begin{theo}\label{arapura} The characteristic variety $V_k(\pi_n(X-D))$ 
is a union of translated subgroups $S_j$ of the group 
${\rm Char} \pi_1(X-D)$
by unitary characters $\rho_j$: $$V_k(\pi_n(X-D))=\bigcup \rho_jS_j$$
\end{theo}

This is an immediate consequence of the interpretation \ref{locsys}  
and the following theorem 
applied to a resolution $\hat X$ 
of non-normal crossings of $D$:  
\begin{theo} (Arapura, \cite{Arapura}) Let $\hat X$
 be a projective manifold such that $H^1(\hat X,\CC)=0$.
Let $\hat D$ be a divisor with normal crossings. Then there
exists a finite number of unitary characters $\rho_j \in {\rm Char} 
\pi_1(\hat X-\hat D)$ and  holomorphic maps $f_j: \hat X-\hat D
\rightarrow T_j$ into complex tori $T_j$
such that the set 
$\Sigma^k(\hat X-\hat D)=\{ \rho \in {\rm Char} \pi_1(\hat X-\hat D)
\vert {\rm dim} H^k(\hat X-\hat D,\rho)  \ge 1  \} $
coincides with $\bigcup \rho_jf^*_iH^1(T_j,\CC^*)$. In particular,
$\Sigma^k$ is a union of translated by unitary characters 
subgroups of ${\rm Char}\pi_1(X-D)$.
\end{theo} 

The components of $\Sigma^1$ can {\it all} be obtained using the maps
$X-D$ onto the curves with negative Euler characteristics 
(cf. \cite{Arapura}). In the case $k>1$, maps onto quasi-projective 
algebraic varieties with abelian fundamental group and 
vanishing $\pi_i$ for $2 \le i \le k-1$ allow one to construct 
components of $V(\pi_k)$ (cf. example \ref{8_4} below).

\section{Review of local theory of isolated non-normal crossings}
\label{review}

Local theory of isolated singularities of 
holomorphic functions provides a beautiful interplay between 
algebraic geometry and topology and in particular the 
topology of (high dimensional) 
links (cf. \cite{milnor}). The main structure is the Milnor fibration 
$\partial B_{\epsilon}-V_f^0 \cap \partial B_{\epsilon} \rightarrow S^1$,
where $V_f^0$ is the zero set of a holomorphic function $f(x_1,..,x_{n+1})$
and $B_{\epsilon}$ is 
a ball of a small radius $\epsilon$ about $\cal O$ (the fibration exist 
even in the non-isolated case). 
If the singularity of $f$ at $\cal O$ 
is isolated, then the fiber $M_f$ of this fibration (the Milnor fiber) is  
homotopy equivalent to a wedge of spheres: $S^n \vee ... \vee S^n$.
Going around the circle, which is the base of Milnor's fibration, 
yields the monodromy: $H^n(M_f) \rightarrow H^n(M_f)$. It  
has as its eigenvalues only the roots of unity $exp(2 \pi i \kappa) \ \
(\kappa  \in \QQ)$. Moreover, there are several ways to pick 
a particular value of the logarithm $\kappa$ of an eigenvalue of the monodromy
so that the corresponding rational number will have 
some geometric significance. One of the ways to do this 
depends on the existence of a Mixed Hodge structure (cf. \cite{Oslo}) on 
$H^n(M_f)$. The value of the logarithm is selected so that 
its integer part is determined by the degree of the component of $Gr^F_*H^n(M_f)$
(graded space associated with the Hodge filtration) on which 
particular eigenvalue of the semi-simple part of the 
monodromy appears.

Some of the data above can be obtained by considering 
the infinite cyclic 
cover of $\partial B_{\epsilon}-V_f^0 \cap \partial B_{\epsilon}$
instead of Milnor fibration.
Such a cover is well-defined since 
$H_1(\partial B_{\epsilon}-V_f^0 \cap \partial B_{\epsilon},\ZZ)=\ZZ$ 
for $n>1$. For example, the universal cyclic cover is diffeomorphic 
to the product $M_f \times \RR$. The monodromy can be identified with 
the deck transformation of the infinite cover. 

With such reformulation, the Milnor theory can be extended 
to the case of germs of isolated non-normal crossings in $\CC^{n+1}$ 
(cf. \cite{innc}), i.e. germs of functions $f_1 \cdot ... \cdot f_r$
such that the intersection points of divisors $f_1=0,...,f_r=0$  
are normal crossings except for the origin $\O$ 
(more general case of germs of complex spaces with isolated 
singularities considered in (cf. \cite{loctop}). The results, 
using infinite covers as a substitute for the Milnor fiber, are parallel to 
the above mentioned results in the isolated singularities case.
Notice, however, that though the theory of Milnor fibers is applicable 
to germs of INNC, much less detailed information can be obtained
since these singularities are not isolated for $n>1$.
For example, the Milnor fiber is not even simply-connected
(cf., below however, where quite a bit of information about 
the Milnor fiber can be obtained as a consequence of 
the present approach).

Let $D$ be a germ of INNC which belongs to a ball $B_{\epsilon}$ 
about $\O$ and which has $r$ irreducible components.
We have the isomorphism $H_1(\partial B_{\epsilon}-D,\ZZ)=\ZZ^r$
and hence the universal abelian cover of $\partial B_{\epsilon}-D$
has $\ZZ^r$ as the covering group. The replacement of the Milnor 
fiber in this abelian situation is the universal abelian 
cover $ \widetilde {\partial B_{\epsilon}-D}$.
Notice that a locally trivial fibration of 
$\partial B_{\epsilon}-D$ over a torus does not exist in general
since typically $\widetilde{\partial B_{\epsilon}-D}$
has the homotopy type of an infinite complex.
We have the following (cf. \cite{innc}):

\begin{theo}\label{INNC} For $n>1$, the fundamental group 
$\pi_1(\partial B_{\epsilon}-D )$ is free abelian.
The universal (abelian) cover $ \widetilde {\partial B_{\epsilon}-D }$
is $(n-1)$-connected. In particular, $H_n(\widetilde {\partial B_{\epsilon}-D},\ZZ)$
is isomorphic to the homotopy group $\pi_n(\partial B_{\epsilon}-D)$.
The latter isomorphism is the isomorphism of $\ZZ[\pi_1(\partial B_{\epsilon}-D)]$-
modules where the module structure on the homology is given by the 
action of $\pi_1(\partial B_{\epsilon}-D)$ on the universal cover via deck 
transformations and the 
action on the homotopy is given by the Whitehead product (cf. \cite{Spanier}) 
\end{theo}
Notice that the case when $D$ is a divisor with normal crossings
is ``a non-singular'' case since the universal cover is contractible.
The simplest example of INNC is given in $\CC^{n+1}$ by the equation 
$l_1 \cdot ...\cdot l_r=0$, where $l_i$ are {\it generic} linear forms
(i.e. a cone over a generic arrangement of 
hyperplanes in $\PP^n$). Since the complement to a generic arrangement
of $r$ hyperplanes in $\PP^n$ has a homotopy type of $n$-skeleton of 
the product of $r-1$-copies of the circle $S^1$ (in minimal cell decomposition
in which one has $r-1 \choose i$ cells of dimension $i$) one can calculate 
the module structure on the $\pi_n$ of such skeleton. Its universal 
cover is obtained by removing the $\ZZ^{r-1}$ orbits of all open faces of 
a dimension greater than $n$ in the unit cube in $\RR^{r-1}$. Hence 
$\pi_n(\partial B_{\epsilon}-D)=H_n(\widetilde{\partial B_{\epsilon}-D},\ZZ)$ 
($\widetilde{\partial B_{\epsilon}-D}$ is the universal cover).
The chain complex of the universal cover of $(S^1)^{r-1}$ can be 
identified with the Koszul complex of the group ring of 
$\ZZ^{r-1}=\ZZ^r/(1,...,1)$ (so that the generators of $\ZZ^{r}$ correspond to the
standard generators of $H_1(\partial B_{\epsilon}-D)$ ). The system of 
parameters of this Koszul complex is $(t_1-1,..,t_r-1)$.
Hence $H_n(\widetilde{\partial B_{\epsilon}-D},\ZZ)=
Ker \Lambda^nR \rightarrow \Lambda^{n-1}R$ 
where $R=\ZZ[t_1,..,t_r]/(t_1 \cdot ...\cdot t_r-1)$. 
As a result, one has the following presentation:
\begin{equation}
\Lambda^{n+1}({\bf Z}[t_1,t_1^{-1},...,t_r,t_r^{-1}]/(t_1...,t_r-1)^r)
\rightarrow \Lambda^n({\bf Z}[t_1,t_1^{-1},...,t_r,t_r^{-1}]/(t_1...,t_r-1)^r)
\rightarrow 
\end{equation}
$$\pi_n({\bf C}^{n+1}-\bigcup D_i) \rightarrow 0 $$
In particular, the support of the $\pi_n$ is the subgroup 
$t_1 \cdot ... \cdot t_r=1$.

We summarize the similarities between the case of hypersurfaces
with isolated singularities and INNC 
in the following table (with \ref{INNC} justifying the first three
rows):

$$\matrix{{\bf Isolated \ singularities} & {\bf INNC} \cr
                         & \cr
                 {\rm Milnor } & {\rm Infinite } \cr
                     {\rm fiber} & {\rm abelian \ cover}\cr 
    & \cr  
{\rm Homology \ of } & {\pi_n } \cr
 {\rm Milnor \ fiber} & \cr
  & \cr 
 {\rm  Monodromy}  & \pi_1-{\rm module \  structure} \cr 
 & {\rm on} \  \pi_n \cr  
 &  \cr
{\rm Eigenvalues} & { \rm Characteristic} \cr 
{\rm of \  monodromy}
& {\rm varieties \ of \  \pi_n} \cr 
& \cr 
{\rm Monodromy}   & {\rm Translated} \cr
{\rm theorem} &  {\rm subgroup \ Property} \cr
&  \cr
{\rm Multiplier} & {\rm Multivariable \ Ideals \ of} \cr
  {\rm Ideals} & {\rm Quasiadjucntion} \cr 
       & \cr
{\rm Spectrum} & {\rm Polytopes \ of} \cr
  & {\rm quasiadjunction} \cr }
$$

\bigskip
\bigskip
In the case of isolated singularities one has the isomorphism:
$\pi_n(\partial B_{\epsilon}-D)=H_n(\widetilde{\partial B_{\epsilon}-D})$
as $\ZZ[t,t^{-1}]$-modules, where the module structure on the 
right is given by the monodromy action. In particular, it is 
a torsion module and its support is a subset of $\Char \ZZ=\CC^*$ 
consisting of the eigenvalues of the monodromy of Milnor 
fibration. Monodromy theorem (\cite{milnor}) is equivalent to 
the assertion that eigenvalues are the torsion points 
of $\CC^*$. A generalization of this is the following:

\begin{conj}\label{transsubgr} The support of $\pi_n(\partial B_{\epsilon}-D)$
is a union of translated subgroups of $\Char \pi_1(\partial B_{\epsilon}-D)$
by points of finite order.
\end{conj}

\ref{transsubgr} is a local analog of the result \ref{arapura}
in the quasi-projective case. 
Now let us describe a partial result in the direction 
of \ref{transsubgr} describing some components of a characteristic 
variety which satisfy \ref{transsubgr}, and which also will explain 
last two rows in the above table. 
  
As already was mentioned, the cohomology of the Milnor fiber $H^n(M_f,\CC)$
of an isolated singularity support a Mixed Hodge structure (cf. \cite{Oslo}).
The monodromy splits into the product of the semi-simple and the unipotent 
part. The semi-simple part leaves the Hodge filtration invariant. 
The latter allows one to split the eigenvalues into groups corresponding to 
the components of $Gr^FH^n(M_f,\CC)$, depending on the graded piece
on which the eigenvalue appears. As a consequence,
 one can assign a rational number
to each eigenvalue, i.e., its logarithm so that its integer part 
is determined by the group to which the eigenvalue belongs 
(we refer to \cite{Oslo} for the 
exact description). In other words, we obtain a lift of the support of 
the homotopy group of the Milnor fiber into the universal cover of the 
subgroup of unitary characters of $\ZZ$ (the eigenvalues of the monondromy
having a finite order are unitary).

In the abelian (local) case, we have the following. Let us consider
the universal cover of the subgroup $\Char^u(\pi_1(\partial B_{\epsilon}-D))$
of unitary characters. It is isomorphic to $\RR^r$ and one can 
take the unit cube as the fundamental domain of the covering group (i.e. 
$\ZZ^r$). We assign an element in the fundamental domain
to a unitary character $\chi$ having finite 
order using the following 
interpretation of the unitary characters from 
$V_k(\pi_n(\partial B_{\epsilon}-D))$ (cf. \cite{innc} Prop. 4.5).  

\begin{prop}\label{localbranch} Let $G=\oplus_{i=1}^{i=r} \ZZ/m_i\ZZ$ 
be a finite quotient of $\pi_1(\partial B_{\epsilon}-D)$
and let $\chi \in \Char(\pi_1(\partial B_{\epsilon}-D))$ which lifts 
to a character of $G$. Then the link $X_{m_1,..,m_r}$ 
of the isolated complete intersection
singularity: 
\begin{equation}\label{completeint}z_1^{m_1}=f_1(x_1,...,x_{n+1}),..., 
z_r^{m_r}=f_{r}(x_1,...,x_{n+1})
\end{equation} is a $n-1$-connected $2n+1$-manifold, 
which is a cover of $\partial B_{\epsilon}$ branched over INNC $D$.
The condition: $\chi \in V_k(\pi_n(\partial B_{\epsilon}-D))$ and 
$\chi$ is essential (cf. \ref{essential}) is equivalent to 
  $$k=\dim \{v \in H_n(X_{m_1,...,m_r}) \vert gv=\chi(g)v \forall g \in G \}$$ 
\end{prop}

Note that the covering map $X_{m_1,..,m_r} \rightarrow \partial B_{\epsilon}$
is just a projection 
$(z_1,...,z_r,x_1,..,x_{n+1}) \rightarrow (x_1,...,x_{n+1})$.
Next, we shall use the Mixed Hodge structure on the cohomology 
of the link (\ref{completeint}) (cf. \cite{arcatamhs}). 
The Hodge filtration $$F^0H^n(X_{m_1,...m_r}) 
 \supset ... \supset F^nH^n(X_{m_1,...m_r}) \supset 0$$
is preserved by the group $G$. The logarithms of characters which appear 
on the subspace $F^nH^n(X_{m_1,...m_r})$ (i.e., 
the vectors $\log\chi=(\xi_1,..,\xi_r)$ with $0 \le \xi_i < 1, \forall i$ such 
that $exp(2\pi i \xi_1),...,\exp(2\pi i \xi_r)$ is
a character $\chi$ of $H_1(\partial B_{\epsilon}-D)$ in coordinates given 
by the generators $H_1(\partial B_{\epsilon}-D)$)
form a polytope in the sense of the following

\begin{dfn} A polytope in the unit cube 
${\cal U}=\{ {\bf x}=(x_1,...,x_n) 
\vert 0 \le x_i \le 1 \forall i \}$ is a subset of $\cal U$ 
formed by the solutions of a system of inequalities 
${\bf a_k} \cdot {\bf x} \le c_k$ for some constants $c_k$ 
(resp. vectors $\bf a_k$) 
such that ${\bf a_k}=(a_k^1,...,a_k^i,...,a_k^n), 0 \le a_k^i \in \QQ$ 
and $0 \le c_k \in \QQ, 
\forall i,k$. A face of a polytope $\cal P$ is a subset of its boundary 
$\partial \cal P$ which has the form $\partial {\cal P} \cap H$ for 
a hyperplane $H$ different from one of $2n$ hyperplanes $x_i=0,1$.
\end{dfn}

We have the following: 

\begin{theo}\label{polytopestheo} To each germ of INNC $D$ and
$l, 0 \le l \le n$ corresponds a collection ${\bf P}_l$ of 
polytopes $\P_{k,l} \in {\bf P}_l$  
such that a vector $\log\chi \in \QQ^r$ in unit cube
belongs to one of the polytopes $\P_{k,l}$ if and only if 
$\dim \{v \in F^l/F^{l+1}H^n \vert gv=\chi(g)v \}=k$. 
In particular, $V_k(\pi_n(\partial B_{\epsilon}-D))=\bigcup_k exp \P_{k,l}$
where $\RR^r \rightarrow \Char^uH_1(\partial B_{\epsilon}-D)$ 
is the exponential map.
\end{theo}

In the cyclic case, each of $\P_{l,k}$ is a rational number $\xi$
such that $exp 2 \pi i \xi$ is an eigenvalue of the monodromy 
having a multiplicity $k$ which appears on $F^l/F^{l+1}H^n(M_f)$, i.e.
is an element of the spectrum having a multiplicity $k$ (in the 
case $l=n$ one obtains the constant of quasi-adjunction 
from  \cite{arcatame}). In the case $n=1$, i.e., the case of 
reducible plane curves, we have the polytopes of quasi-adjunction 
studied in \cite{alexhodg}. In particular, these polytopes 
are related to the multi-variable log-canonical thresholds and 
multiplier ideals (cf. remark 2.6 and section 4.2 respec. in 
\cite{alexhodg}). Similar relations exist in the case 
of INNC discussed here. 

In the case of isolated singularities, there are very explicit and 
beautiful calculations of the eigenvalues of the monodromy and
spectrum of singularities. We would like to pose the following
problem: 

\begin{prob} Calculate the characteristic varieties of INNC with 
$\CC^*$-actions and in the case when $D_i$ are generic for 
their Newton polytopes. What are the polytopes described in the 
theorem \ref{polytopestheo}?
\end{prob}

This should be a generalization of the case, 
discussed above, of the cone over a generic arrangement 
and the example in \cite{innc} of the cone over a divisor with normal 
crossings in $\PP^n$.

\section{Homology of abelian covers}\label{covers}

In this section, we return to the global case of divisors with ample components
having only isolated non-normal crossings.

\subsection{Topology of unbranched covers} 

The characteristic varieties $\Char_i(\pi_n(X-D))$
contain information about both branched and unbranched abelian covers.

\begin{lem}\label{fincovlocsys}
 Let $G$ be a finite abelian quotient of $\pi_1(X-D)$
and let $U_G$  be  corresponding 
unbranched covers
of $X-D$. 
Let $\chi \in Char(\pi_1(X-D))$ be a pull back of 
a character of $G$  (we shall considered it as a character of the 
latter). Let $H_n(U_G)_{\chi}=\{ v \in H_n(U_G) \vert
g \cdot v=\chi (g)v ( g \in G)\}$ Then 
$H_n(X-D,\L_{\chi})=H_n(U_G)_{\chi}$. In particular,  
$\chi \in {\rm Char} G \subset 
{\rm Char}_i(\pi_n)$ if and only if $H_n(U_G)_{\chi} \ge i$.
\end{lem}

A proof can be obtained, for example, from the exact sequence
of low degree non-vanishing terms in the 
spectral sequence of the action of the group 
$K={\rm Ker} \pi_1({X-D}) \rightarrow G$
on the universal cover $\widetilde {X-D}$ (for which 
we have $\widetilde {X-D}/K=U_G$):
$$H_p(K,H_q(\widetilde{X-D})) \Rightarrow H_{p+q}(U_G)$$
(cf. \cite{cartaneilenberg}).
This is a spectral sequence of $\CC[\pi_1(X-D)]$-modules
where the $\CC[\pi_1(X-D)]$-module structure on
$\CC[G]$ module $H_{p+q}(U_G)$ comes via surjection:
$\CC[\pi_1(X-D)] \rightarrow  \CC[G]$. The localization 
of this spectral sequence at 
a point $\chi$ of ${\rm Char}G \subset {\rm Char}\pi_1(X-D)$
yields the claim using \ref{locsys}, since the localization of 
$H_n(U_G)$ at $\chi$ has the same $\chi$-eigenspace 
as $H_n(U_G)$.

Now, let us consider the effect of adding (ample) components to $D$.
\begin{lem}\label{addingcomponents}
Let $D'$ an ample divisor such that $D \cup D'$ is a divisor
with isolated non-normal crossings. Then the homomorphism of 
$\pi_1(X-D)$ modules: $\pi_i(X-D\cup D') 
\rightarrow \pi_i(X-D)$ is surjective for $1 \le i \le \dim X-1$.
In particular, if one considers $\Spec\CC[\pi_1(X-D)]$ as a subset 
in $\Spec\CC[\pi_1(X-D\cup D')]$, then the intersection of 
$V_k(\pi_n(X-D \cup D'))$ with  $\Spec\CC[\pi_1(X-D)]$ contains 
$V_k(\pi_n(X-D))$.
\end{lem}
{\it Sketch of the proof}. Let $T(D')$ be a small neighborhood 
of $D'$ in $X$. Then by the Lefschetz theorem, $\pi(T(D')-D'\cap D)$ 
surjects onto $\pi_i(X-D \cup D')$. On the other hand, 
this map can be factored through $\pi_i(X-D)$ which yields the claim.
\smallskip
\par \noindent Lemma \ref{addingcomponents} suggests the following definition:
\begin{dfn}\label{essential} The components of $V_k(X-\bar D)$ 
where $\bar D$ is a union of a collection of $D_i$'s forming
$D$ and which are 
considered as subsets in $\Spec \CC[\pi_1(X-D)]$ called 
the non-essential components of $V_k(X-D)$. The remaining components 
are called essential.

 A character $\chi$ is called essential 
if $\chi(\gamma) \ne 1$ for each element $\gamma \in \pi_1(X-D)$ 
which is a boundary of a small 2-disk transversal to one of 
irreducible components of $D$.
\end{dfn}
We shall see in the next section that only essential characters 
contribute to the homology of branched covers.

\subsection{Hodge theory of branched covers.} The relationship between 
the homology of branched and unbranched covers is more subtle in the present 
case than in the case of plane curves considered in 
\cite{abcov} and the local case of the section \ref{review}. 
One of the reasons is that there is no prefer non singular model for 
the abelian global case. Only the birational type of branched cover
is an invariant of $X-D$, and hence the Betti numbers of branched covers 
depend upon compactification of the unbranched cover. However, the Hodge 
numbers $h^{i,0}$ are birational invariants (in the case ${\rm dim}X=2$, 
they determine the relevant part of homology of branched cover 
completely due to the relation $b_1=2h^{1,0}$) and one can expect a relation
between the Hodge numbers $h^{i,0}$ and the homology of 
unbranched covers. 

Recall that the cohomology of unitary local systems supports a 
mixed Hodge structure (cf. \cite{timm}). We shall denote
$h^{p,q,k}(\L)={\rm dim} Gr^p_FGr^q_{\bar F}Gr^W_{p+q}(H^k(\L))$ 
the dimension  
of the corresponding Hodge space.
In the case of a rank one  
local system having a finite order, one has the following 
counterpart of \ref{fincovlocsys}:

\begin{theo}\label{branchedcover}
Let, as in \ref{fincovlocsys}, 
$\chi \in Char(\pi_1(X-D))$ be a character of 
a finite quotient $G$ of $\pi_1(X-D)$. Let $\bar U_G$ 
be a $G$-equivariant non-singular compactification of $U_G$ and
let $H^{p,q}(\bar U_G)_{\chi}$ be the $\chi$-eigenspace of $G$ 
acting on $H^{p,q}(\bar U_G)$. 
Then 
$$ h^{n,0,n}(\L_{\chi})=h^{n,0}(\bar U_G)_{\chi}$$
\end{theo}

\noindent {\it Sketch of the Proof.} The functoriality of 
the Hodge structure on cohomology of local systems yields
that the isomorphism in \ref{locsys} is compatible with the 
Hodge structure:
$h^{n,0,n}(X-D,\L_{\chi})=h^{n,0,n}(U_G)_{\chi}$ where 
in the RHS are the Hodge numbers of the Deligne's MHS 
on the cohomology of non-singular quasi-projective manifold 
(cf. \cite{DeligneII}). Let $E=\bar U_G-U_G$, which we assume is a  divisor 
with normal crossings. In the exact sequence of MHS:
$H^n(\bar U_G,U_G) \rightarrow H^n (\bar U_G) \rightarrow H^n(U_G)$,
which splits into corresponding sequences of $\chi$-eigenspaces,
the image of right homomorphism is $W_nH^n(U_G)$ (cf. \cite{DeligneII} 3.2.17).
This result is a consequence of the 
identity: ${\rm Ker} H^n (\bar U_G) \rightarrow H^n(U_G) \cap H^{n,0}=0$
To see the latter, notice that using the duality
$H^{n+2}(E) \times H^n(\bar U_G,U_G) \rightarrow \CC(-n-1)$
($\CC(-k)$ is the Hodge-Tate) 
we obtain $h^{n,0,0}(\bar U_G,U_G)=h^{n+1,1,n+2}(E)$.
On the other hand, for each smooth component $E_i$ of $E$
one has $h^{i,j,n+2} \ne 0$ only when $0 \le i,j \le n$ and 
the Mayer Vietoris sequence of MHS yields the same conclusion
for $E$. Hence $H^n (\bar U_G) \rightarrow H^n(U_G)$ is injective
on $H^{n,0}$ and the result follows.   

\section{Conjecture and results on the structure of 
characteristic varieties.}\label{conjectures}

In this section we shall continue to study the 
situation which was discussed in section \ref{action} and consider 
complements to divisor $D$ with isolated non-normal crossings
on projective manifolds $X$ ($\dim X=n+1$). Our goal 
is to calculate components of $V_i(\pi_n(X-D))$. 
The procedure described below is a generalization of 
the one outlined in \cite{abcov}.

First, we shall associate with each divisor in $Pic(X) \otimes \RR$
a subspace in the universal cover of the group of unitary characters 
of $\pi_1(X-C)$.
Let us assume that 
$H_1(X-D)=\ZZ^r$ (i.e., to avoid notational complications, assume it 
is torsion free) and let $l_i (i=1,...,r)$ be a system of generators
of $H_1(X-D)$. The boundary operator in (\ref{homologycomplement}) 
associates to each $D_j$ an element 
$\delta_j=\sum a_{i,j}l_i (a_{i,j} \in \ZZ$. 
In particular, elements of $\Char H_1(U,\ZZ)$ define
the elements in $\Hom (H^{2n}(D),\CC^*)=H_{2n}(D,\CC^*)$ 
whose image in $H_{2n}(X,\CC^*)=H^2(X,\CC^*)$ 
belongs to $Pic(X)\otimes \CC^*$ (which can be identified 
with the subgroup of the latter since $\pi_1(X)=0$). Similarly 
the elements in the universal cover of the maximal compact 
subgroup of $\Char H_1(X-D,\ZZ)$, which we identify with 
$\Hom(H_1(X-D,\ZZ),\RR)=H^1(X-D,\RR)$
define the elements in the
universal cover of the maximal compact subgroup of 
$H^{2n}(D,\CC^*)$ and  hence in the maximal compact 
subgroup of $Pic (X) \otimes \CC^*$. The last universal 
cover will be identified with $Pic(X) \otimes \RR$.
The lift into universal cover of the 
elements of finite order in $\Char H_1(X-D,\ZZ)$ 
yields the elements in $Pic(X) \otimes \QQ$. 
Preimage $L_D$ of an element $D \in Pic(X) \otimes \RR$ in 
$H^1(U,\RR)$ is an affine subspace in the universal cover of the 
space of characters. For example, if $X=\PP^{n+1}$, $D$ is union 
of hypersurfaces of degrees $1,d_1, ...,d_r$ (i.e., $H_1(\PP^{n+1}-D)=
\ZZ^{r}$), one has $Pic(X)=\ZZ$ and the preimage of $\O(l) \in Pic(X)$ 
is the hyperplane in $H^1(X-D,\RR)$ given by $x_1d_1+...+x_rd_r=l$ 
where $x_i$'s are the coordinates corresponding to the basis of $H_1(X-D,\ZZ)$ 
given by the cycles dual to $D_i$'s.
 
Now, with each $S \in \S$, we associate a polytope in 
the unit cube in $\RR^r$ as follows.
For any $S \in \S$, one has the map 
$H_1(B_{\epsilon}-D) \rightarrow H_1(X-D)$ and hence
the map $\Char H_1(X-D) \rightarrow \Char H_1(B_{\epsilon}-D)$.
The latter lifts to the map of universal covers: $\RR^r \rightarrow \RR^s$
where $s$ is the number of components of $D$ containing $S$.
This can be described in coordinates as follows.
A vector $\Xi: (\kappa_1,..,\kappa_r) \in \QQ^r (0 \le \kappa_i < 1)$
for any collection $(j_1,...,j_s)$ determines the vector:
$\Xi^{j_1,..,j_s}=
(\{\sum a_{i,j_1}\kappa_i\},....,\{\sum a_{i,j_1}\kappa_i\} ) \in \QQ^r$
($\{ \dot \}$ is the fractional part of a rational number).
For each $S \in \S$, we consider subsets $\P_S^{gl}$ consisting 
of vectors $\Xi=(\kappa_1,..,\kappa_r)$ such that 
$\Xi^{j_1,..,j_s} \in \P_S$ where $D_{j_1},.,D_{j_s}$ are the 
components of $D$ passing through $S \in \S$ and $\P_S \in \QQ^s$ is 
a face of a polytope of quasi-adjunction of INNC formed by 
$D_{j_1},..,D_{j_s}$.

\begin{dfn} Let $\S \subset X$ be the collection 
of non-normal crossings of the 
divisor $D$.
Global polytope of quasi-adjunction 
corresponding to $\S$ is $\bigcap_{S \in \S} \P_S^{gl}$.
A global face of quasi-adjunction is a face of a global polytope of
quasiadjunction.
A divisor $\D=\sum \alpha_i D_i \in Pic X, \alpha_i \in \QQ$ 
is called contributing if the corresponding subset $L_{\D}$ 
of the elements of 
the universal cover (cf. above) contains  
a global face of quasiadjunction $\F$ and 
$H^1(\A_{\F} \otimes \Omega^{n+1}_X \otimes \D) \ne 0$.
A global face of quasiadjunction $\F$ is contributing if there 
is a contributing divisor such that the corresponding 
subspace $L_{\D}$ contains $\F$.
\end{dfn}
 
\begin{conj}Zariski closure of $exp(\F)\subset \Char H_1(U,\ZZ)$
is a component of characteristic variety $V_k$ where
$k=\dim H^1(\A_{\F} \otimes \Omega^{n+1}_X \otimes \D)$.
\end{conj}

I don't know if such components are {\it all} essential components 
if the characteristic variety.

The supporting evidence is the following. This conjecture is shown 
in \cite{abcov} in the case of curves and in the case $X=\PP^{n+1}$
and $D$ is a hypersurface with isolated singularities in 
\cite{Maninvolume}. Both of these 
results can be generalized as follows (the proof, based 
on the methods used in these two papers will appear elsewhere).

\begin{theo} Let $D \subset \PP^{n+1}$ be a union of hypersurfaces
$D_0,D_1,...,D_r$ of degrees $1,d_1,...,d_r$ respectively, 
which is a divisor with isolated 
non-normal crossings. Let $\F$ be a face of global polytope of 
quasi-adjunction, i.e. 
a face of an intersection of polytopes of quasi-adjunction corresponding 
to a collection $\S$ of 
 non-normal crossings of $D$. Let $d_1x_1+...+d_rx_r=l$ be
a hyperplane containing the face of quasiadjunction $\F$.
If $H^1(\A_{\F}\otimes \O(l-3)) =k$, then the Zariski closure of 
$exp (\F) \subset \Char H_1(\PP^{n+1}-D)$ is a component of 
$V_k(\pi_n(X-D))$.
 
\end{theo}

A consequence of the conjecture is the corollary.

\begin{conj} Let $\chi \in \Char \pi_1(X-D)$ be a character of a finite
quotient of $G$ of the fundamental group. 
Let, as in lemma \ref{branchedcover}, $\bar U_G$ be a $G$-equivariant 
compactification of the unbranched cover of $X-D$ with the Galois 
group $G$ and let $h^{n,0}_{\chi}(\bar U_G)$ be the $\chi$-eigenspace 
of the $G$ acting on $H^{n,0}(U_G)$. Then  $h^{n,0}_{\chi}(\bar U_G)=0$
unless the lift of $\chi$ belongs to a contributing global 
face of quasi-adjunction $\F$ and in which case one has:
 $$h^{n,0}_{\chi}(\bar U_G)=\dim H^1(\A_{\F} \otimes \Omega^n_X \otimes
\O(D))$$ 
where $D \in Pic(X)$ is the divisor corresponding to the lift of 
character $\chi$.
\end{conj}

In the case when there are bundles $\L_i$ such that $\L_i^{n_i}={\cal O}(D_i)$ and the 
cover corresponds to the group $\ZZ_{n_1} \times ... \times \ZZ_{n_r}$,
using the arguments similar to those used in \cite{Vaquie},
one obtains (in agreement with the conjecture)
the following: the eigenspace of the action of $G$ corresponding 
to the dimension of the 
eigenspace corresponding to $(e^{2 \pi i {{p_1} \over {n_1}}},...,
e^{2 \pi i {{p_r} \over {n_r}}})$ 
is $\dim H^1(\A_{\F} \otimes \Omega^{n+1}_X \otimes \L_1^{p_1} \otimes ...
\otimes \L_r^{p_r})$.
In the case when $r=1$ the condition that $({{p_1} \over {n_1}},...,
{{p_r} \over {n_r}})$ belongs to the face of quasi-adjunction 
becomes the condition that ${p \over n}$ is an element of the spectrum of 
one of the singularities of the divisor $D$ and one obtains the result
from \cite{Vaquie}.

\section{Examples}\label{examples}

\subsection{Local examples} 
\begin{exam} Germs of curves.
\end{exam}
In the case of curves, the support 
of $H_1(\widetilde {\partial B_{\epsilon}-D},\CC)$ is the zero set 
of the Alexander polynomial. There are extensive calculations 
of this invariant using knot-theoretical methods (cf. \cite{EN}). 
Hodge decomposition is considered in \cite{alexhodg}.
For example, for the singularity $x^r-y^r=0$, the characteristic variety 
is $t_1 \cdot \cdot \cdot t_r=1$ (cf. the calculation for the cone 
over the generic arrangement in section \ref{review}). The 
faces of the polytopes of quasi-adjunction are the hyperplanes 
$x_1+...+x_r=l, \ (l=1,..,r-2)$.

\begin{exam} Cones.
\end{exam}
A generalization of the example of 
arrangements of hyperplanes considered in section \ref{review} 
is given by a union of non singular hypersurfaces in ${\bf P}^n$ 
which form a divisor with normal crossings (cf. \cite{innc}). 
If the degrees of hypersurfaces
are $d_1,..,d_r$ respectively then $V_1={\rm Supp}(\pi_n(\CC^{n+1}-D) \otimes \CC)$ 
is given by $t_1^{d_1} \cdot ... \cdot t_r^{d_r}-1=0$.

\subsection{Global examples} 

\begin{exam} Plane curves
\end{exam} 

We refer to \cite{abcov} for examples of characteristic varieties for 
pencils quadrics (Ceva arrangement of four lines) and pencils of cubics (arrangement 
of nine lines dual to inflection points of a non-singular cubic and 
the arrangement of 12 lines containing its inflection points).
Papers \cite{LY1}, \cite{LY2} and \cite{CS} describe a combinatorial 
method to detect components of characteristic variety
and in \cite{jose} a generalization to arrangements 
of rational curves is considered. Papers \cite{artal} and \cite{bandman} 
contain applications of 
characteristic varieties to geometric problems. 

\begin{exam}\label{8_4} Arrangement in $\PP^3$ with
 isolated non-normal crossings for which 
$\pi_2$ of the complement which support has non-trivial essential 
components.
\end{exam}

Consider the  
arrangement $D_{8,4}$ of hyperplanes in ${\bf P}^3$ 
which is an (8,4) configuration (cf. \cite{Gonz}). 
It includes a plane containing 4 generic points $Q_1,...,Q_4$,
six generic planes $H_{i,j}$ each passing through the line $Q_iQ_j$ and 
also the plane containing the four coplanar (by Desargue theorem) 
points $H_{i,j} \cap H_{j,k} \cap H_{i,k}$.
Recall (cf. \cite{Gonz}) that this configuration contains 
eight planes and eight points such that every plane contains 
four points and every point belongs to exactly four planes.
Denoting eight points by $1,2,3,4,1',2',3',4'$ and eight 
planes by ${\bf 1,2,3,4,1',2',3',4'}$ the incidence relation is 
given by the diagram: 
$$\matrix{1 & 1' & & {\bf 1} & {\bf 1'} \cr 
         2 & 2' & & {\bf 2} & {\bf 2'} \cr
              3 & 3' & & {\bf 3} & {\bf 3'} \cr
       4 & 4' & & {\bf 4}  &{\bf 4'} \cr }$$

\noindent where the plane in position $(i,j)$ 
contains all points in row $i$ and column $j$ except for the 
point in position $(i,j)$. 

This arrangement of eight hyperplanes has only 
isolated non-normal crossings.
From \ref{h1}, we infer that $H_1(\PP^3-D_{8,4})=\ZZ^7$.
Moreover we have the rational 
map: 
$$\Pi: \ \  \PP^3 \rightarrow \PP(H^0(\PP^3,{\cal I}(2)))^*=\PP^2$$
where $\cal I$ is the ideal sheaf of the collection of eight points 
in $\PP^3$ forming this configuration. The indeterminancy points 
are the eight points of configuration. In order to calculate the $\Pi$-image
of the hyperplanes of the arrangement, notice that the points in the target of 
the map correspond to the pencils of quadrics in the web, 
the image of a point is the pencil of quadric in the web containing 
this point and the lines correspond to quadrics in $ H^0(\PP^3,{\cal I}(2))$
i.e. are the collections of pencils containing a quadric. 
In particular, the image of a point $P$ in 
a hyperplane $H \in D_{8,4}$ is a pencil of quadrics from  
$H^0(\PP^3,{\cal I}(2))$ containing $P$. This pencil contains the  
quadric among the four quadrics containing $P$, mentioned earlier. 
Hence the image of 
$P$ belongs to the union $L$ of four lines in $\PP(H^0(\PP^3,{\cal I}(2))^*$
corresponding to above four quadrics. Therefore we have 
a regular map: $\Pi: \PP^3-D_{8,4} \rightarrow \PP^2-L$.

Let us calculate the cohomology of local systems $\Pi^*(\L)$, where
$\L$ is a local system on $\PP^2-L$. We have the Leray spectral sequence:
$H^p(\PP^2-L,R\Pi_*^q(\Pi^*\L)) \Rightarrow H^{p+q}(\PP^3-D_{8,4},\Pi^*\L)$.
In it $R\Pi_*^q(\L)=\L \otimes R^q_*(\CC)$ i.e. 
$H^0(\PP^2-L,\L \otimes R\Pi^1_*(\CC))=0$ for all but finitely 
many local systems and $\Pi_*(\CC)=\CC,
R\Pi_*^2(\CC)=0$ (the latter since the fiber of $\Pi$ is not compact).
Hence this spectral sequence degenerates.
Moreover, since $L$ is a generic arrangement and hence 
$\pi_1(\PP^2-L)$ is abelian, for any non-trivial local system on $\PP^2-L$, 
the first cohomology is vanishing and one has 
$H^1(\PP^2-L,R\Pi_*^1(\Pi^*\L))=0$. Hence $H^2(\PP^3-D_{8,4},\Pi^*(\L))=
H^2(\PP^2-L,\CC)=\CC^3$. 

The above calculation shows that, since 
$\pi_2(\PP^2-L) \otimes \CC=\CC[t_1,t,_2,t_3,t_4]/(1-t_1t_2t_3t_4)$,
the support of the homotopy group: 
$\pi_2(\PP^3-D_{8,4})\otimes \CC$ has a 3-dimensional 
component. Projections from each of eight vertices of this configuration 
yield linear maps of $\PP^3-D_{8,4}$ onto the complement in $\PP^2$ 
to  four lines in a general position and hence a 3-dimensional 
component. We obtain hence in ${\rm Spec}\CC[H_1(\PP^3-D_{8,3})]$
nine 3-dimensional components of the support of $\pi_2(\PP^3-D_{8,3})
\otimes \CC$. 

\section{Betti and Hodge realizations of multi-variable motivic zeta function}

The purpose of this section is to relate the motivic zeta function 
of Denef and Loeser in the case of local INNC  
to the invariants considered in section (\ref{review}).

Recall that, to a smooth variety $X$ over $\CC$ and 
$r$ holomorphic functions $f_i: X \rightarrow \CC$, one associates a 
multivariable motivic zeta-function $Z_{f_1,...,f_r}(T_1,..,T_r)$
which is a formal
series in 
$\M_{X_0 \times \GG_m^r}[[T_1,...,T_r]]$. Here, as in 
 \cite{loeser0006}, $X_0=\bigcap_i f_i^{-1}(0)$, $\GG_m$ is the 
multiplicative group of the field $\CC$ and 
for a variety  
$S$ the ring $\M_S$ is obtained from the Grothendieck group 
$K_0(Var_S)$
of varieties over $S$ 
by inverting the class $\LL$ of $\AA^1_k \times S \in K_0(Var_S)$.
More precisely, denote $\L(X)$ (resp. $\L_n(X)$) the arc
space of $X$ (resp. arc space mod $n$) whose points are 
the maps $\Spec \CC[[t]] \rightarrow X$ (resp.
 $\Spec \CC[[t]]/(t^{n+1}) \rightarrow X$). Let  
$${\cal X}_{n_1,...,n_r}=\{\phi \in \L_n(X), n=\sum n_j 
\vert {\rm ord}_t \phi^*(f_j)=n_j
 \ \ j=1,...,r  \}$$
and $ac(f)=(ac(f_1),...,ac(f_r)):{\cal X}_{n_1,..,m_r} \rightarrow \GG_m^r$
assigns to an arc in ${\cal X}_{n_1,..,n_r}$
the vector which $j$-th component is the coefficient of $t^{n_j}$ in 
$\phi^*(f_j)$. Together with $\pi_0: {\cal X}_{n_1,..,n_r} \rightarrow X$ 
which 
assignes to an arc the image in $X$ of its closed point 
$\Spec \CC \rightarrow \Spec \CC[[t]]$, this makes ${\cal X}_{n_1,..,n_r}$ 
into $\GG_m^r \times X_0$-manifold.
Then 
\begin{equation}\label{motiviczeta}
 Z_{f_1,...,f_r}(T_1,..,T_r)=\sum_{n_1,..,n_r, n_i \in \bf N}
[{\cal X}_{n_1,...,n_r}/X_0 \times \GG_m^r]\LL^{(d\sum n_i)}T_1^{n_1} \cdot
\cdot \cdot T_r^{n_r}  
\end{equation}
One has the canonical 
maps (resp. Betti and Hodge realizations): 
 $e_{top}: K_0(Var_{\CC}) \rightarrow \ZZ$ and 
$e_h: K_0(Var_{\CC}) \rightarrow \ZZ[u,v]$ induced by the maps 
assigning to a variety $V$ its topological euler characteristic
and the E-function $\sum_i (-1)^i\dim Gr^p_FGr^W_{p+q}H^i(V)u^pv^q$
(both $F$ and $W$ filtration are coming from Deligne's Mixed Hodge 
structure on $V$). 
We also will use the equivariant refinement  
of $e_{top}$ and $e_h$ defined for  
$V \in Var_{\CC}$ supporting an action 
of a finite group $G$ 
via biholomorphic transformations. For $\chi \in \Char G$,
those refinements pick the corresponding eigenspaces: 
\begin{equation} e_{top,\chi}(V)=\sum (-1)^i \dim H^i(V)_{\chi} \ \ \ 
{\rm and} \ \
e^m_{h,\chi}(V)=\sum_i (-1)^i Gr^m_FH^i(V)_{\chi}
\end{equation}
The function (\ref{motiviczeta}) can be expressed 
in terms of a resolution of singularities of $f_1,,,,f_r$ as follows
(cf. \cite{loeser0006}).
Let $Y \rightarrow X$ be a resolution of singularities of $D$, 
i.e. the union of the exceptional set $\bigcup_{i \in J} E_i$ and the 
proper preimage of an INNC $D$ is a normal crossings divisor. 
For $I \subset J$, let 
$E_I^{\circ}=\cap_{i \in I}E_i-\cup_{j \in J-I} E_j$, 
$a_{i,k}$ (resp. $c_k$) is the order along the exceptional 
component of the pull-back on $Y$ of function $f_i$ (resp. the 
order of the pull back of the differential 
$dx_1 \wedge... \wedge dx_{n+1}$). 
Let $U_i$ be the complement to the zero section of the normal 
bundle to $E_i$ in $Y$, and $U_I$ is the fiber product of 
$U_i \vert_{E_I^{\circ}}$ over $E_I$.
Then:

\begin{equation}\label{zetaresolution}
Z_{f_1,..f_r}(T_1,..,T_r)=\sum_{I \subset J}[U_J/\GG^r_m \times X_0]
\Pi_{i \in I}
{{\LL^{-c_i-1}T_1^{a_{i,1}} \cdot \cdot \cdot T_r^{a_{i,r}}} \over
{1-\LL^{-c_i-1}T_1^{a_{i,1}} \cdot \cdot \cdot T_r^{a_{i,r}}}}
\end{equation}

We have the following:

\begin{theo} Betti realization of $Z_{f_1,...,f_r}(T_1,...,T_r)$ 
determines the essential components of 
the characteristic variety $V_1$. More precisely, for 
an essential $\chi$: 
\begin{equation}\label{limit} 
V_1=-e_{top,\chi}
{\lim_{T_i \rightarrow \infty}}Z_{f_1,...,f_r}(T_1,...,T_r)
\end{equation}
\end{theo}
{\it Proof}. One can deduce this 
from C.Sabbah's results in \cite{sabbah} similarly to \cite{guibert}
since, due to the vanishing theorem \ref{INNC}, the 
multi-variable zeta function 
studied in \cite{sabbah} determines the support of the 
$\pi_1(B_{\epsilon}-D)$ module $\pi_n(B_{\epsilon}-D)$.

\bigskip In the cyclic case, the Hodge realization of the motivic
zeta function is equivalent to the spectrum (cf. \cite{loeser0006}).
At least in the case of curves, one has the Hodge version 
in the abelian case as well (as was suggested in \cite{loesercastling}):

\begin{theo} For $n=1$, the Hodge realization of 
(\ref{motiviczeta}) determines the polytopes of quasiadjunction.
\end{theo}
{\it Proof}. Let $X_{m_1,...,m_r}$ be the link of an abelian cover 
$\V_{m_1.,,.,m_r}$ given by the equations (\ref{completeint})
with $n=1$. A resolution of this complete 
intersection singularity in the category of spaces with quotient singularities
(in the case of surfaces with ADE singularities) can be obtained as the 
normalization $\widetilde {V_{m_1,..,m_r}}$ 
of $\V_{m_1,..,m+r} \times_{B_{\epsilon}} Y_D$, where
$Y_D \rightarrow B_{\epsilon}$ is an embedded resolution of the 
singularities of $D$. The exceptional locus $\tilde E$ of 
the resolution of (\ref{completeint}) supports the action 
of the group $G=\ZZ_{m_1} \times ...\times \ZZ_{m_r}$.
We have the following sequence of MHS (cf. \cite{arcatamhs}): 
$0 \rightarrow H^1_E(\widetilde {V_{m_1,..,m_r}}) \rightarrow H^1(E) 
\rightarrow H^1( X_{m_1,...,m_r}) \rightarrow 0$, which in the case $n=1$
yields the equivariant isomorphism $H^1( X_{m_1,...,m_r})=H^1(\tilde E)$ 
of MHSs. Since the MHS on $H^2(\tilde E)$ is pure,
we have:  
\begin{equation}\label{lastone}
\dim F^1H^1(X_{m_1,..,m_r})_{\chi}=\dim Gr^1_FH^1(L)_{\chi}=
e^{1}_{h,\chi}(\tilde E) 
\end{equation}
The latter is determined by the Hodge realization of (\ref{zetaresolution}),
since the pull-back of $[U_i]$ via the map $\M_{\GG^r_m \times X_0}
 \rightarrow \M_{\GG^r_m \times X_0}$ corresponding to the 
map ${\GG^r_m \times X_0} \rightarrow {\GG^r_m \times X_0}$ given by 
 $z_i=u_i^{m_i}$ 
is equivalent to the unbranched cover of 
$\partial B_{\epsilon}-D$, which is 
preimage of $\GG_m^r \subset \CC^r$ for the projection of (\ref{completeint})
onto the space of $z$-coordinates. 
In particular, it determines the class of the exceptional set $\tilde E_i$ 
in $\M_{\CC}$.
It follows from (\ref{lastone}) that
$\dim F^1H^1(X_{m_1,..,m_r})_{\chi} \ge 1$ iff 
$e^{1}_{h,\chi}(\tilde E) \ge 1$ QED.


\begin{thebibliography}{99}

\bibitem{Abh} S.Abhyankar, Tame Covering of fundamental groups of algebraic 
varieties I,II, Amer. Journ. of Math. vol. 81, 1959 p.46-94, 
vo. 82, 1960.

\bibitem{Arapura} D.Arapura, Geometry of cohomology support loci for local systems. I.  J. Algebraic Geom.  6  (1997),  no. 3, 563--597. 

\bibitem{artal} Artal, E.; Carmona, J.; Cogolludo, J. I.; Tokunaga, Hiro-O Sextics with singular points in special position.  J. Knot Theory Ramifications  10  (2001),  no. 4, 547--578. 

\bibitem{symplectic} D. Auroux, S. K. Donaldson, L. Katzarkov, M. Yotov,
Fundamental groups of complements of plane curves and symplectic invariants,
(math.GT/0203183).

\bibitem{bandman} Bandman T, Libgober, A. Counting rational maps onto surfaces
and fundamental groups. Preprint. 2004.

\bibitem{Brown} K. Brown, Cohomology of groups. Graduate Texts in Mathematics, 87. Springer-Verlag, New York-Berlin, 1982

\bibitem{cartaneilenberg} H.Cartan, S.Eilenberg, Homological algebra. 
Princeton University Press, Princeton, N. J., 1956

\bibitem{denis} Ch\'eniot, D.; Libgober, A. Zariski-van Kampen theorem 
for higher-homotopy groups.  J. Inst. Math. Jussieu  2  (2003),  no. 4, 495--527.
 
\bibitem{jose} J.I.Cogolludo, Topological Invariants of the Complements to 
Rational Arrangements, Thesis, University of Illinois at Chicago, 1999.

\bibitem{CS} D.Cohen, A.Suciu,  Characteristic varieties of arrangements.  
Math. Proc. Cambridge Philos. Soc. 127 (1999), no. 1, 33--53.

\bibitem{DeligneII} P.Deligne, Th\'eorie de Hodge. II. (French)  Inst. Hautes Etudes Sci. Publ. Math. No. 40 (1971), 5--57.

\bibitem{loeser0006} Denef, J. Loeser, F. 
 Geometry on arc spaces of algebraic varieties.  
European Congress of Mathematics, Vol. I (Barcelona, 2000),  327--348, Progr. Math., 201, Birkh\"user, Basel, 2001. 
Geometry of arc spaces of 
algebraic varieties.


\bibitem{DP} Dimca, A. Papadima, S.
 Hypersurface Complements, Milnor Fibers and 
Higher Homotopy Groups of Arrangements.
 Ann. of Math. (2)  158  (2003),  no. 2, 473--507.

\bibitem{loctop} Dimca, A. Libgober, A. Local topology of reducible divisors,
(math.AG/0303215) 

\bibitem{EN} Eisenbud, D. Neumann, W. Three-dimensional link theory and invariants 
of plane curve singularities. Annals of Mathematics Studies, 110. Princeton University Press, Princeton, NJ, 1985. 

\bibitem{Fulton} Fulton, W. On the topology of algebraic varieties.  Algebraic geometry, Bowdoin, 1985 (Brunswick, Maine, 1985),  15--46, Proc. Sympos. Pure Math., 46, Part 1, Amer. Math. Soc., Providence, RI, 1987.

\bibitem{Gonz} Gonzalez-Dorrego, M. (16.6)-configurations and Geometry 
of Kummer surfaces in $\PP^3$. Memoirs AMS, vol. 107, 152, (1994).

\bibitem{GM} Goresky, M, MacPherson, R.
 Stratified Morse theory. Ergebnisse der Mathematik und ihrer Grenzgebiete (3) [Results in Mathematics and Related Areas (3)], 14. Springer-Verlag, Berlin, 1988. 

\bibitem{guibert} Guibert, G. Espaces d'arcs et invariants d'Alexander. 
Comment. Math. Helv.  77  (2002),  no. 4, 783--820. 


\bibitem{arcatame} Libgober, A. Alexander invariants of plane 
algebraic curves.  Singularities, Part 2 (Arcata, Calif., 1981), 
 135--143, Proc. Sympos. Pure Math., 40, Amer. Math. Soc., Providence, RI, 1983. 


\bibitem{arcatame2} Libgober, A. Fundamental groups of the complements 
to plane singular curves.  
Algebraic geometry, Bowdoin, 1985 (Brunswick, Maine, 1985), 
 29--45, Proc. Sympos. Pure Math., 
46, Part 2, Amer. Math. Soc., Providence, RI, 1987. 

\bibitem{Annals} Libgober, A. Homotopy groups of 
the complements to singular hypersurfaces. II.  
Ann. of Math. (2)  139  (1994),  no. 1, 117--144. 

\bibitem{Maninvolume} Libgober, A. Position of 
singularities of hypersurfaces and the topology of their complements. 
Algebraic geometry, 5.  J. Math. Sci.  82  (1996),  no. 1, 3194--3210. 


\bibitem{abcov} Libgober, A. Characteristic varieties of algebraic curves. 
 Applications of algebraic geometry to coding theory, physics and computation 
(Eilat, 2001),  215--254, NATO Sci. Ser. II Math. Phys. Chem., 36,
 Kluwer Acad. Publ., Dordrecht, 2001. 

\bibitem{alexhodg}  Libgober, A. Hodge decomposition of Alexander invariants.  Manuscripta Math.  107  (2002),  no. 2, 251--269. 


\bibitem{innc} Libgober, A. Isolated non normal crossings Proceedings
of the workshop on Singularities, San Paolo, 2002. M.Ruis and T.Gaffney Editors, 
(math.AG/0211264)

\bibitem{LY1} Libgober, A, Yuzvinsky, S, Cohomology of the Orlik-Solomon 
algebras and local systems.  Compositio Math.  121  (2000),  no. 3, 337--361. 

\bibitem{LY2}  Libgober, A. Yuzvinsky, S. Cohomology of local systems.  
Arrangements---Tokyo 1998,  169--184, Adv. Stud. Pure Math., 27, 
Kinokuniya, Tokyo, 2000


\bibitem{loesercastling} Loeser, F. 
Motivic zeta functions for prehomogeneous vector spaces and 
castling transformations.  Nagoya Math. J.  171  (2003), 85--105

\bibitem{milnor} Milnor, J, Singular points of complex hypersurfaces.
 Annals of Mathematics Studies, No. 61 Princeton University Press, 
Princeton, N.J.; University of Tokyo Press, Tokyo 1968 

\bibitem{Nori}  Nori, M, Zariski's conjecture and related problems.
Ann. Sci. \"Ecole Norm. Sup. (4) 16 (1983), no. 2, 305--344.

\bibitem{SP} Papadima, S., Suciu, A. 
Higher homotopy groups of complements of complex hyperplane arrangements.  Adv. Math.  165  (2002),  no. 1, 71--100.

\bibitem{sabbah}  Sabbah, C, Modules d'Alexander et $\D$-modules. 
Duke Math. J.  60  (1990),  no. 3, 729--814.

\bibitem{Spanier} Spanier, E. Algebraic Topology, McGraw Hill Book Company, 
1966.

\bibitem{Oslo}  Steenbrink, J. H. M. Mixed Hodge structure on the 
vanishing cohomology.  Real and complex singularities 
(Proc. Ninth Nordic Summer School/NAVF Sympos. Math., Oslo, 1976), 
 pp. 525--563. Sijthoff and Noordhoff, Alphen aan den Rijn, 1977. 


\bibitem{arcatamhs}  Steenbrink, J. H. M. Mixed Hodge structures 
associated with isolated singularities.  Singularities, Part 2 
(Arcata, Calif., 1981),  513--536, Proc. Sympos. Pure Math., 40, 
Amer. Math. Soc., Providence, RI, 1983. 

\bibitem{mina} Teicher, M. Braid groups, algebraic surfaces and fundamental groups of complements of branch curves.  Algebraic geometry---Santa Cruz 1995,  127--150, Proc. Sympos. Pure Math., 62, Part 1, Amer. Math. Soc., Providence, RI, 1997.

\bibitem{timm} Timmerscheidt, K. Mixed Hodge theory for unitary local systems.  J. Reine Angew. Math.  379  (1987), 152--171.

\bibitem{Vaquie} Vaqui\'e, Michel Irr\"egularit\"e des rev\"etements 
cycliques. Singularities (Lille, 1991),  383--419, 
London Math. Soc. Lecture Note Ser., 201, Cambridge Univ. Press, Cambridge, 1994.

\end{thebibliography}
\end{document}